\documentclass{llncs}
\usepackage{makeidx}  % allows for indexgeneration

\usepackage{epsfig}
\usepackage{amssymb}

\begin{document}

\title{Many Touchings Force Many Crossings}
%\title{Hamiltonian Mechanics unter besonderer Ber\"ucksichtigung der
%h\"ohreren Lehranstalten}
%
\titlerunning{Many Touchings Force Many Crossings}
%\titlerunning{Hamiltonian Mechanics}  % abbreviated title (for running head)

\author{J\'anos Pach\inst{1,2} and G\'eza T\'oth\inst{2}}

%\author{Ivar Ekeland\inst{1} \and Roger Temam\inst{2}
%Jeffrey Dean \and David Grove \and Craig Chambers \and Kim~B.~Bruce \and
%Elsa Bertino}
%

\authorrunning{J. Pach, G. T\'oth}

%\authorrunning{Ivar Ekeland et al.} % abbreviated author list (for running head)
%Craig Chambers, Kim B. Bruce, and Elisa Bertino}
%

\institute{
\'Ecole Polytechnique F\'ed\'erale de Lausanne, St.  8, Lausanne 1015, Switzerland\\
\email{pach@cims.nyu.edu}
\and
R\'enyi Institute, Hungarian Academy of Sciences 1364 Budapest, POB 127, Hungary\\
\email{ geza@renyi.hu}
}

%\institute{Princeton University, Princeton NJ 08544, USA,\\
%\email{I.Ekeland@princeton.edu},\\ WWW home page:
%\texttt{http://users/\homedir iekeland/web/welcome.html}
%\and
%Universit\'{e} de Paris-Sud,
%Laboratoire d'Analyse Num\'{e}rique, B\^{a}timent 425,\\
%F-91405 Orsay Cedex, France}

\maketitle              % typeset the title of the contribution

\begin{abstract}
Given $n$ continuous open curves in the plane, we say that a pair
is {\em touching} if they have only one interior point in common
and at this point the first curve does not get from one side of
the second curve to its other side. Otherwise, if the two curves
intersect, they are said to form a {\em crossing} pair.
Let $t$ and $c$ denote the number of touching pairs and crossing pairs,
respectively. We prove that $c \ge {1\over 10^5}{t^2\over n^2}$,
provided that $t\ge 10n$. Apart from the values of the constants, this 
result is best possible.
\end{abstract}

\section{Introduction}

In the context of the theory of topological graphs and graph drawing,
many interesting questions have been raised concerning the adjacency
structure of a family of curves in the plane or in another surface
\cite{FP10}. In particular, during the past four decades, various
important properties of string graphs (i.e., intersection graphs
of curves in the plane) have been discovered, and the study of
different crossing numbers of graphs and their relations to one another
has become a vast area of research. A useful tool in these investigations
is the so-called crossing lemma of Ajtai, Chv\'atal, Newborn,
Szemer\'edi and Leighton
\cite{ACNS}, \cite{L}. It states the following: Given a graph of $n$ vertices
and
$e>4n$ edges, no matter how we draw it in the plane by not necessarily
straight-line
edges, there are at least constant times $e^3/n^2$ crossing pairs of edges.
\smallskip

This lemma has inspired a number of results establishing the existence
of many crossing subconfigurations of a given type in sufficiently
rich geometric or topological structures \cite{D98}, \cite{Sh03}, 
\cite{SoT01}, \cite{GNT00}.

\smallskip

In this note, we will be concerned with families of curves in the plane.
By a {\em curve}, we mean a non-selfintersecting continuous arc in the plane,
that is, a homeomorphic image of the open interval $(0,1)$. 
Two curves are said to
{\em touch} each other if they have precisely one interior point in common and at this
point the first curve does not pass from one side of the second curve to the 
other.
Any other pair of curves with nonempty intersection is called {\em crossing}.
A family of curves is in {\em general position} if any two of them intersect
in a finite number of points and no three pass through the same point.

\smallskip

Let $n$ be even, $t$ be a multiple of $n$, and suppose that 
$n\le t<{n^2\over 4}$.
Consider a collection $A$ of $n-{2t\over n}>{n\over 2}$ pairwise disjoint
curves, and another collection $B$ of ${2t\over n}$ curves such that

(i) $A\cup B$ is in general position,

(ii) each element of $B$ touches precisely ${n\over 2}$ elements of $A$, and

(iii) no two elements of $B$ touch each other.

\noindent The family $A\cup B$ consists of $n$ curves such that the number of
touching pairs among them is $t$. The only pairs of curves that may cross each
other belong to $B$. Thus, the number of crossing pairs is at most ${2t/n
  \choose 2}\le {2t^2\over n^2}$. See Figure 1.

\begin{figure}[h]
 \begin{center}
  \includegraphics[scale=0.45]{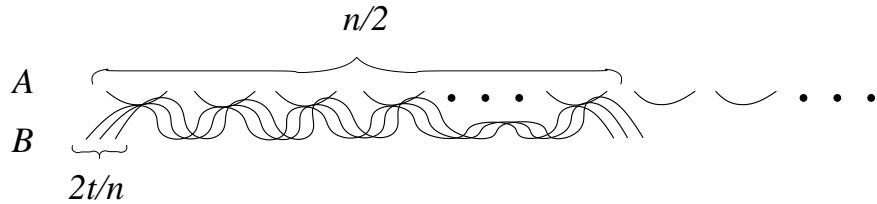}
 \end{center}
\caption{A set of $n$ curves with $t$ touching pairs and at most 
${2t^2\over n^2}$ crossing pairs.}
\end{figure}

The aim of the present note is to prove that this construction is optimal up
to a constant factor, that is, any family of $n$ curves and $t$ touchings has
at least constant times ${t^2\over n^2}$ crossing pairs.

\medskip

\noindent{\bf Theorem.} {\em 
Consider a family of $n$ curves in general position
in the plane which determines $t$ touching pairs and $c$ crossing pairs.

If $t\ge 10n$, then we have $c \ge {1\over 10^5}{t^2\over n^2}$.
This bound is best possible up to a constant factor.}
\medskip

We make no attempt to optimize the constants in the theorem.
\smallskip

Pach, Rubin, and Tardos
\cite{PRT16} established a similar relationship between $t$,
the number of touching pairs, and $C$, the number of crossing {\em points}
between the curves. They proved that $C \ge t(\log\log (t/n))^{\delta}$,
for an absolute constant $\delta>0$. Obviously, we have $C\ge c$.
There is an arrangement of $n$ red curves and $n$ blue curves in the
plane such that every red curve touches every blue curve,
and the total number of crossing points is $C=\Theta(n^2\log n)$;
cf.~\cite{FFPP10}.
Of course, the number of crossing pairs, $c$, can never exceed ${n\choose 2}$.

Between $n$ arbitrary curves, the number of touchings $t$ can be as large as
$(\frac34 + o(1)){n\choose 2}$; cf. \cite{PT06}. However, if we restrict our
attention to algebraic plane curves of bounded degree, then we have
$t=O(n^{3/2})$, where the constant hidden in the notation depends on the
degree \cite{ESZ16}.

\section{Proof of Theorem}

We start with an easy observation.
\medskip

\noindent {\bf Lemma.} {\em Given a family of $n\ge 3$ curves in general
  position
in the plane, no two of which cross, the number of touchings, $t$, cannot
exceed $3n-6$.}

\medskip

\noindent {\em Proof.}
Pick a different point on each curve. Whenever two curves touch each
other at a point $p$, connect them by an edge (arc) passing through $p$.
In the resulting drawing, any two edges that do not share an endpoint
are represented by disjoint arcs. 
According to the Hanani-Tutte theorem~\cite{Tu70},
this means that the underlying graph is planar, so that its number of edges,
$t$, satisfies  $t\le 3n-6$. $\Box$

\medskip

\noindent  {\em Proof of Theorem.}
We proceed by induction on $n$. For $n\le 20$, the statement is
void. Suppose that $n>20$ and that the statement has already
been proved for all values smaller than $n$.

We distinguish two cases.
\medskip

\noindent{\sc CASE A:} $t \le 10n^{3/2}$.
\smallskip

In this case, we want to establish the stronger statement
$$c\ge {1\over 10^4}{t^2\over n^2}.$$
By the assumption, we have
\begin{equation}\label{eq0}
{1\over 10^4}{t^2\over n^2}\le {n\over 100}.
\end{equation}

Let  $G_t$ (resp., $G_c$) denote the {\em touching graph}
(resp., {\em crossing graph}) associated with the curves.
That is, the vertices of both graphs correspond to the curves,
and two vertices are connected by an edge if and only if the
corresponding curves are touching (resp., crossing).

\smallskip

%If $G_t$ has a vertex of degree smaller than $t/n$, delete
%it. The resulting touching graph has $n-1$ vertices and $t'$ edges,
%where
%${t'\over (n-1)}\ge {t\over n}$. Thus, by the induction hypothesis,
%among the remaining $n-1$ curves, there are  at least
%${1\over 10^3}{t'^2\over (n-1)^2} \ge {1\over 10^3}{t^2\over n^2}$
%crossing pairs, so we are done. Therefore, in the sequel,
%we can and will assume that
%the degree of each vertex in $G_t$ satisfies
%\begin{equation}\label{eq1}
%{\rm deg}_{G_t}(v)\ge{t\over n}.
%\end{equation}
%\smallskip

\vskip 0.3cm

\begin{figure}[h]
 \begin{center}
  \includegraphics[scale=0.4]{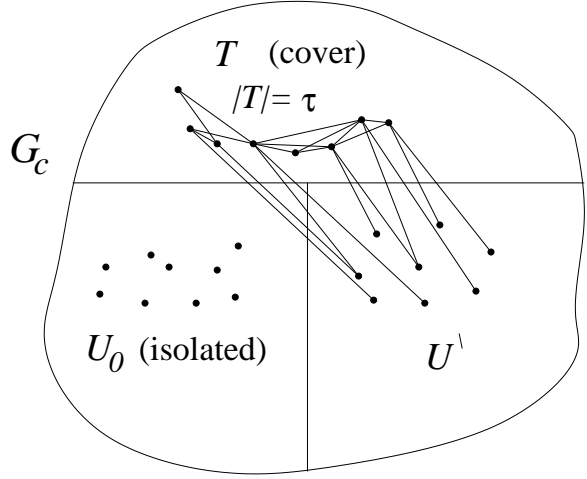}
 \end{center}
\caption{Graph $G_c$.}
\end{figure}

Let $T$ be a minimal vertex cover in  $G_c$, that is, a
smallest set of vertices of $G_c$ such that
every edge of $G_c$ has at least one endpoint in $T$. Let
$\tau=|T|$.
Let $U$ denote the complement of $T$.
Obviously, $U$ is an {\em independent set} in $G_c$.
According to the Lemma,  the number of edges in 
$G_t[U]$, the touching graph induced by $U$, satisfies
\begin{equation}\label{eq1.5}
|E(G_t[U])|< 3|U|\le 3n.
\end{equation}

By the minimality of $T$,  $G_c$ has at least $|T|=\tau$ edges. That is, we have
$c \ge \tau$, so
we are done if
$\tau\ge {1\over 10^4}{t^2\over n^2}$.
\smallskip

From now on, we can and shall assume that
$\tau < {1\over 10^4}{t^2\over n^2}.$
By (\ref{eq0}), we have ${1\over 10^4}{t^2\over n^2}\le {n\over 100}$. 
Hence, $|T|\le {n\over 100}$ and
\begin{equation}\label{eq2.5}
|U|=n-|T| \ge {99n\over 100}.
\end{equation}

Let $U'\subseteq U$ denote the set of all vertices in $U$ that are not isolated
in the graph $G_c$. By the definition of $T$, all neighbors 
of a vertex $v\in U$ in $G_c$  belong to $T$.
If $|U'|\ge {1\over 10^4}{t^2\over n^2}$,
then we are done, because $c\ge |U'|.$
\smallskip

Therefore, we can assume that
\begin{equation}\label{eq3}
|U'|< {1\over 10^4}{t^2\over n^2}\le {n\over 100},
\end{equation}
where the second inequality follows again by (\ref{eq0}).

Letting $U_0=U\setminus U'$,
by (\ref{eq2.5}) and (\ref{eq3}) we obtain $|U_0|=|U|-|U'|\ge {98n\over 100}$.
Clearly, all vertices in $U_0$ are isolated in $G_c$.

Suppose that $G_t[T\cup U']$ has at least ${t\over 10}$ edges.
Consider the set of curves $T\cup U'$. We have
$n_0=|T\cup U'|\le {2n\over 100}$ and, the number of touchings,
$t_0=|E(G_t[T\cup U'])|\ge {t\over 10}$. Therefore,
by the induction hypothesis, for the number of crossings we have
$c_0=|E(G_c[T\cup U'])|\ge 
{1\over 10^5}{t_0^2\over n_0^2}\ge {1\over 10^4}{t^2\over n^2}$
and we are done.
Hence, we assume in the sequel that
$G_t[T\cup U']$ has fewer than ${t\over 10}$ edges.

Consequently, for
the number of edges in $G_t$ running between $T$ and $U_0$, we have
\begin{equation}\label{eq4}
|E(G_t[T,U_0])|\ge t-|E(G_t[T\cup U'])|-|E(G_t[U_0\cup U'])|\ge 
t-{t\over 10}-3n>{t\over 2}.
\end{equation}
Here we used the assumption that $t\ge 10n$.

Let $\chi=\chi(G_c[T])$ denote the chromatic number of $G_c[T]$.
In any coloring of a graph with the smallest possible number of colors,
there is at least one edge between any two color classes.
Hence,  $G_c[T]$ has at least 
${\chi\choose 2}\ge {1\over 10^4}{t^2\over n^2}$ edges, 
and we are done, provided that
$\chi>{1\over 70}\cdot{t\over n}$.
\smallskip

Thus, we can suppose that
\begin{equation}\label{eq5}
\chi=\chi(G_c[T])\le {1\over 70}\cdot{t\over n}.
\end{equation}

Consider a coloring of  $G_c[T]$ with $\chi$ colors, and denote
the color classes by $I_1, I_2, \ldots , I_{\chi}$.
Obviously, for every $j$, $I_j\cup U_0$ is an independent set in $G_c$.
Therefore, by the Lemma, $G_t[I_j\cup U_0]$ has at most $3n$ edges.
Summing up for all $j$ and taking (\ref{eq5}) into account, we obtain
$$|E(G_t[T,U_0])|\le \sum_{j=1}^{\chi}|E(G_t[I_j\cup U_0])|
\le{1\over 70}\cdot{t\over n}3n\le {t\over 20},$$
contradicting (\ref{eq4}). This completes the proof in CASE A.

\medskip

\noindent{\sc CASE B:} $t\ge 10n^{3/2}$.
\smallskip
Set  $p={10n^3\over t^2}\le {1 \over 10}$.
Select each curve independently with probability $p$.
Let ${\bf n'}$, ${\bf t'}$, and ${\bf c'}$ denote the number of selected curves,
the number of touching pairs, and the number of crossing pairs between them,
respectively.
Clearly,
\begin{equation}\label{varhatoertek}
E[{\bf n'}]=pn, \;\; E[{\bf t'}]=p^2t, \;\; E[{\bf c'}]=p^2c.
\end{equation}

The number of selected curves, ${\bf n'}$, has binomial distribution, therefore,
\begin{equation}\label{n'}
{\rm Prob}[|{\bf n'}-pn|>{1\over 4}pn]<{1\over 3}.
\end{equation}

By Markov's inequality,
\begin{equation}\label{c'}
{\rm Prob}[{\bf c'}>3p^2c]<{1\over 3}.
\end{equation}

Consider the touching graph $G_t$. Let $d_1, \ldots, d_n$ denote the degrees
of the vertices of $G_t$, and let
$e_1, \ldots , e_t$ denote its edges, listed in any order.
We say that an edge $e_i$ is {\em selected} (or belongs to the random sample)
if both of its endpoints were selected.
Let $X_i$ be the {\em indicator variable} for $e_i$,
that is,
\[
    X_i=\left\{
                \begin{array}{ll}
                  1\;\;\;\mbox{ if $e_i$ was selected,}\\
                  0\;\;\;\mbox{ otherwise.}\\
                \end{array}
              \right.
  \]

We have $E[X_i]=p^2$. Let ${\bf t'}=\sum_{i=1}^tX_i$.
It follows by straightforward computation that for every $i$,
$${\rm var}[X_i]=E[(X_i-E[X_i])^2]=p^2-p^4,$$
If $e_i$ and $e_j$ have a common endpoint for some $i\neq j$, then
$${\rm cov}[X_i, X_j]=E[X_iX_j]-E[X_i]E[X_j]=p^3-p^4.$$
If $e_i$ and $e_j$ do not have a common vertex, then $X_i$ and $X_j$ 
are independent random variables and ${\rm cov}[X_i, X_j]=0$.
Therefore, we obtain
$$\sigma^2={\rm var}[{\bf t'}]
  =\sum_{i=1}^t{\rm var}[X_i]+\sum_{1\le i\neq j\le t}{\rm cov}[X_i, X_j]$$
%$$=\sum_{i=1}^t{\rm var}[X_i]+\sum_{1\le i\neq j\le t}{\rm cov}[X_i, X_j]$$
$$=(p^2-p^4)t+(p^3-p^4)\sum_{i=1}^nd_i(d_i-1)<p^2t+2p^3nt.$$
From here, we get
$\sigma<\sqrt{p^2t}+\sqrt{2p^3nt}<p^2t=E[{\bf t'}]$.
By Chebyshev's inequality,
$${\rm Prob}[|{\bf t'}-p^2t|\ge \lambda\sigma]\le {1\over \lambda^2}.$$
Setting $\lambda={1\over 4}$,
\begin{equation}\label{t'}
{\rm Prob}[|{\bf t'}-p^2t|\ge {p^2t\over 4}]\le {1\over 4^2}<{1\over 3}.
\end{equation}

It follows from (\ref{n'}),  (\ref{c'}), and  (\ref{t'})
that,
with positive probability, we have
\begin{equation}\label{meg0}
|{\bf n'}-pn|\le{1\over 4}pn,\;\;\;\; {\bf c'}\le 
3p^2c,\;\;\;\;|{\bf t'}-p^2t|\le{1\over 4}p^2t.
\end{equation}
\smallskip

Consider a fixed selection of $n'$ curves with $t'$ touching pairs and $c'$
crossing pairs for which the above three inequalities are satisfied.
Then we have
$$t'\ge {3\over 4}p^2t={300\over 4}\cdot{n^6\over t^3},$$
$$n'\le {5\over 4}pn={50\over 4}\cdot{n^4\over t^2},$$
and, hence,
\begin{equation}\label{meg1}
t'\ge {6n^2\over t}n'\ge 10n'.
\end{equation}
On the other hand,
$$t'\le {5\over 4}p^2t={500\over 4}\cdot{n^6\over t^3},$$
$$n'\ge {3\over 4}pn={30\over 4}\cdot{n^4\over t^2},$$
so that
\begin{equation}\label{meg2}
10(n')^{3/2}\ge 10\cdot{30^{3/2}\over 4^{3/2}}\cdot{n^6\over t^3}>t'.
\end{equation}

According to (\ref{meg1}) and (\ref{meg2}), the selected family meets 
the requirements of the Theorem in CASE A.
Thus, we can apply the Theorem in this case to obtain that
$c'\ge {1\over 10^4}{t'^2\over n'^2}$.
In view of (\ref{meg0}), we have
$$3p^2c\ge c',\;\;\;\; t'\ge{3\over 4}p^2t,\;\;\;\; n'\le {5\over 4}pn.$$
Thus,
$$3p^2c\ge c'\ge{1\over 10^4}{t'^2\over n'^2}\ge
{1\over 10^4}{(3p^2t/4)^2\over (5pn/4)^2}
={1\over 10^4}\left({3\over 5}\right)^2{p^2t^2\over n^2}.$$
Comparing the left-hand side and the right-hand side, we conclude that
$$c\ge{1\over 10^5}{t^2\over n^2},$$
as required.
This completes the proof of the Theorem. $\Box$

\bigskip

\noindent{\bf Acknowledgment.} The work of J\'anos Pach was partially 
supported by Swiss National Science Foundation Grants
200021-165977 and 200020-162884. G\'eza T\'oth's work was partially supported
by the Hungarian National Research, Development and Innovation Office, NKFIH,
Grant K-111827.

\end{document}